\documentclass[%
 aip,
 amsmath,amssymb,
 reprint,%
]{revtex4-1}

\usepackage{graphicx}
\usepackage{dcolumn}
\usepackage{bm}

\usepackage[utf8]{inputenc}
\usepackage[T1]{fontenc}
\usepackage{mathptmx}
\usepackage{etoolbox}
\usepackage{subfigure}
\usepackage{graphicx}
\usepackage{dcolumn}
\usepackage{bm}
\usepackage{enumitem}
\usepackage{tabularx} 
\usepackage{enumitem} 

\usepackage{xcolor}
\definecolor{bg}{rgb}{0.975,0.95,1.0}

\definecolor{mygreen}{rgb}{0.1,0.6,0.1}

\usepackage{listings}
\definecolor{codegreen}{rgb}{0,0.6,0}
\definecolor{codegray}{rgb}{0.5,0.5,0.5}
\definecolor{codepurple}{rgb}{0.58,0,0.82}
\definecolor{backcolour}{rgb}{0.975,0.95,1.0}

\lstdefinestyle{mystyle}{
    backgroundcolor=\color{backcolour},
    commentstyle=\color{codegreen},
    keywordstyle=\color{magenta},
    numberstyle=\tiny\color{codegray},
    stringstyle=\color{codepurple},
    basicstyle=\ttfamily\small,
    breakatwhitespace=false,
    breaklines=true,
    captionpos=b,
    keepspaces=true,
    showspaces=false,
    showstringspaces=false,
    showtabs=false,
    tabsize=4,
    numbers=left,
    numbersep=5pt,
}

\makeatletter
\def\@email#1#2{%
 \endgroup
 \patchcmd{\titleblock@produce}
  {\frontmatter@RRAPformat}
  {\frontmatter@RRAPformat{\produce@RRAP{*#1\href{mailto:#2}{#2}}}\frontmatter@RRAPformat}
  {}{}
}%
\makeatother
\begin{document} 

\preprint{AIP/123-QED}

\title[Basins of attraction]{
Effortless estimation of basins of attraction
}

\author{George Datseris}
\email{george.datseris@mpimet.mpg.de}
\affiliation{Max Planck Institute for Meteorology, 20146 Hamburg, Germany}

\author{Alexandre Wagemakers}
\affiliation{Nonlinear Dynamics, Chaos and Complex Systems Group, Departamento de Física, Universidad Rey Juan Carlos, Móstoles, Madrid, Tulipán s/n, 28933, Spain}

\date{\today}

\begin{abstract}
We present a fully automated method that identifies attractors and their basins of attraction without approximations of the dynamics. The method works by defining a finite state machine on top of the dynamical system flow. The input to the method is a dynamical system evolution rule and a grid that partitions the state space. No prior knowledge of the number, location, or nature of the attractors is required. The method works for arbitrarily-high-dimensional dynamical systems, both discrete and continuous. It also works for stroboscopic maps, Poincar\'e maps, and projections of high-dimensional dynamics to a lower-dimensional space. The method is accompanied by a performant open-source implementation in the DynamicalSystems.jl library. The performance of the method outclasses the naive approach of evolving initial conditions until convergence to an attractor, even when excluding the task of first identifying the attractors from the comparison. We showcase the power of our implementation on several scenarios, including interlaced chaotic attractors, high-dimensional state spaces, fractal basin boundaries, and interlaced attracting periodic orbits, among others. The output of our method can be straightforwardly used to calculate concepts such as basin stability and final state sensitivity.
\end{abstract}

\maketitle

\begin{quotation}
Basins of attraction play a central role in the study of multistable dynamical systems. They contain the information about the sets of initial conditions whose trajectories converge to different asymptotic states. Since the basins of most nonlinear dynamical systems are impossible to study analytically, numerical simulations is the method of choice for the inquiry. The computation of the basins implies matching the trajectory of each chosen initial condition against a collection of known attractors. Our algorithm not only automatically identifies the attractors of a dynamical system but also estimates the basins for a given grid of initial conditions.
\end{quotation}

\section{\label{sec:Introduction}Introduction}

In the state space of a dynamical system, basins of attraction are the set of initial conditions that lead to a particular attractor. If only a single global attractor exists, then every initial condition ends up there. However, the coexistence of several attractors in the state space, known as multistability, has been observed in a large array of different dynamical systems\cite{pisarchik2014control}. The recent advent of the tipping-points analysis~\cite{Feudel2018} has enhanced the interest for this phenomenon. In presence of multistability, it is thus important to map the initial conditions to the attractor they end up at, or in other words, to evaluate the basin of attraction of each attractor.

Estimating the basins has benefits well beyond simply knowing the long-term behavior of each initial condition. For example, they can reveal the existence of chaotic transient dynamics before settling into a non-chaotic attractor~\cite{aguirre2009fractal}. Some basins have fractal boundaries. There, it is important to measure how uncertain we are about the final state of an initial condition. It can be computed via different tools, e.g., the uncertainty dimension of the boundary (also known as final state sensitivity)~\cite{grebogi1983final}, or the basin entropy~\cite{daza2016basin,puy2021test}.

Importantly, the basins of attraction can be used to complement or extend the traditional linear stability analysis of the attractors and unveil potential tipping points in a dynamical system~\cite{menck2013basin, Feudel2018}. For example, the basin stability quantifies the robustness of an attractor relative to a perturbation in a system parameter~\cite{menck2013basin}. Also leveraging the information carried by the basins, the tipping probabilities uncover the influence of a parameter drift on the global dynamics~\cite{kaszas2019tipping}. Notice that all the methods we have outlined so far assume that the basins and the attractors have been estimated correctly beforehand.

There are several approaches to construct an approximation of the basins. A brute force method, consisting in evolving initial conditions for long transient and then comparing the last $N$ points of the trajectory, may work well for fixed point attractors. But for anything else it will fail due to the many practical drawbacks regarding, e.g., the variability of integration time-stepping and sampling of non-fixed-point attractors. An alternative is to compare the Lyapunov spectrum of each orbit~\cite{Freire2008} to classify the attractors. The benefit is a simpler comparison between orbits, but it is at the cost of a precise computation of the Lyapunov exponents. Besides, we cannot be sure of the uniqueness of the spectrum for different orbits. For example two symmetric attractors can posses the same spectrum. A third approach relies on recursive subdivision of the state space with quad trees structures~\cite{saupe1987efficient}, that has been useful in estimating basin boundaries of Julia and Mandelbrot sets. However, due to the memory requirement of the quad tree structure, this method is inefficient when the boundary occupies a large portion of the state space or when the state space is higher dimensional. It is unsuitable for generic dynamical systems.

In this article we solve the problem of the computation of the basins by utilizing the only property of an attractor that is always guaranteed to identify it uniquely: its location in the state space. Our approach is inspired by a method described by Nusse \& Yorke in Ref. \onlinecite[Chap. 7]{nusse1997dynamics}. Our algorithm relies on the Poincar\'e recurrence theorem, which states that a trajectory on an attracting set will sooner or later visit the same regions of the state space. The algorithm first locates the attractors by searching for recurrences on a discretized state space grid. The second step is to match initial conditions with attractors, which can be done efficiently both during and after the attractors have been located and labeled. These tasks are executed by pairing the dynamical system with a finite state machine.

We have implemented the algorithm on top of the DynamicalSystems.jl software library \cite{datseris2018dynamicalsystems}, written entirely in the Julia programming language. The algorithm implementation is user-friendly, requiring $\sim$10 lines of input code (these include actually defining the dynamical system). 

In Sect.~\ref{sec:algorithm} we explain the details and potential drawbacks of the algorithm and in Sect.~\ref{sec:results} we apply it successfully on a wide range of different scenarios, from interlaced chaotic attractors to high-dimensional dynamical systems. In Sect.~\ref{sec:implementation} we showcase the code implementation, its computational aspects, as well as the advantages of making it part of a general purpose library. In Sect.~\ref{sec:conclusions} we summarize and conclude.

\section{Description of the algorithm}
\label{sec:algorithm}
\subsection{Attractor identification via recurrences}
\label{sec:attractor_identification}
Identifying attractors using recurrences is a central part of our algorithm and thus we describe it first here in isolation, before moving on to the main algorithm presentation in Sect.~\ref{sec:fsm}.

A portion of the state space of the dynamical system is discretized in the form of a regular grid of initial conditions. 
An array with the same size as the grid is defined for the storage of the information regarding basins and attractors. Each element of this array will hold the information about a cell that is centered around a single initial condition. This information will be called the \emph{label} \texttt{n} of the cell. The size of the cell is determined by the grid step along each dimension. The other component is an integrator that progresses a state space point forwards in time along the flow of the dynamical system.

\begin{figure}
  \includegraphics[width=\columnwidth]{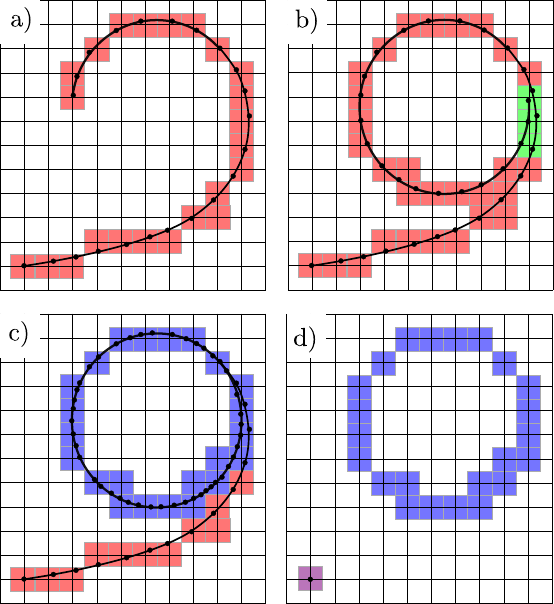}
  \caption{\label{fig:attractor_identification} Attractor identification on the grid. The intersections of the grid correspond to the initial conditions, black dots correspond to states of the dynamical system during integration, and solid lines are a guide the eye (only the black dots are known during the process). The colored areas centered around the intersections are the boxes or cells that will be used for the identification of the attractors and basins. a) As the trajectory evolves the algorithm leaves a mark on each visited cell (red squares). b) When the orbit visits a cell already marked in a), the algorithm begins counting the recurrences (green squares). When the trajectory visits \texttt{mx\_chk\_fnd\_att = 3} consecutive green cells, we consider that we have found a new attractor. c) The algorithm proceeds to locate the attractor correctly. From this moment every visited cell is marked as a part of the attractor (blue squares) and the process goes on until we have visited \texttt{mx\_chk\_loc\_att} blues squares in a row. d) At this point the algorithm erases the marks (red squares) and labels the cell of the initial condition as part of the basin of the attractor (the magenta square).}
\end{figure}

For the task of attractor identification, we track the successive steps of the dynamical system evolution on the grid starting from an initial condition. Each visited cell is labeled \texttt{v} for ``visited'' (red color in Fig.~\ref{fig:attractor_identification}a), and an internal counter registers the sequential events. The trajectory will eventually step \texttt{mx\_chk\_fnd\_att} consecutive times into such \texttt{v}-cells. At this point, we have found an attractor since there are sufficient recurrences on the grid (green color in Fig.~\ref{fig:attractor_identification}b). Note that labels, states, and parameters of the algorithm will be denoted with typewriter characters. Afterwards the algorithm proceeds to locate the rest of the attractor as precisely as possible (the internal counter is reset to 0 here). From this moment on, all visited cells will be marked as containing an attractor point (blue color in Fig.~\ref{fig:attractor_identification}c). 
This encoding goes on until the internal counter reaches \texttt{mx\_chk\_loc\_att}. This ensures that we find as much cells as desired with attractor points. At the end of the process, the algorithm marks the initial condition as part of the basin of the found attractor (magenta color Fig.~\ref{fig:attractor_identification}d). Finally, it discards the labelling \texttt{v} on all other cells visited by the transient of the trajectory.

\subsection{The finite state machine}
\label{sec:fsm}
\begin{figure}[t]
  \includegraphics[width=\columnwidth]{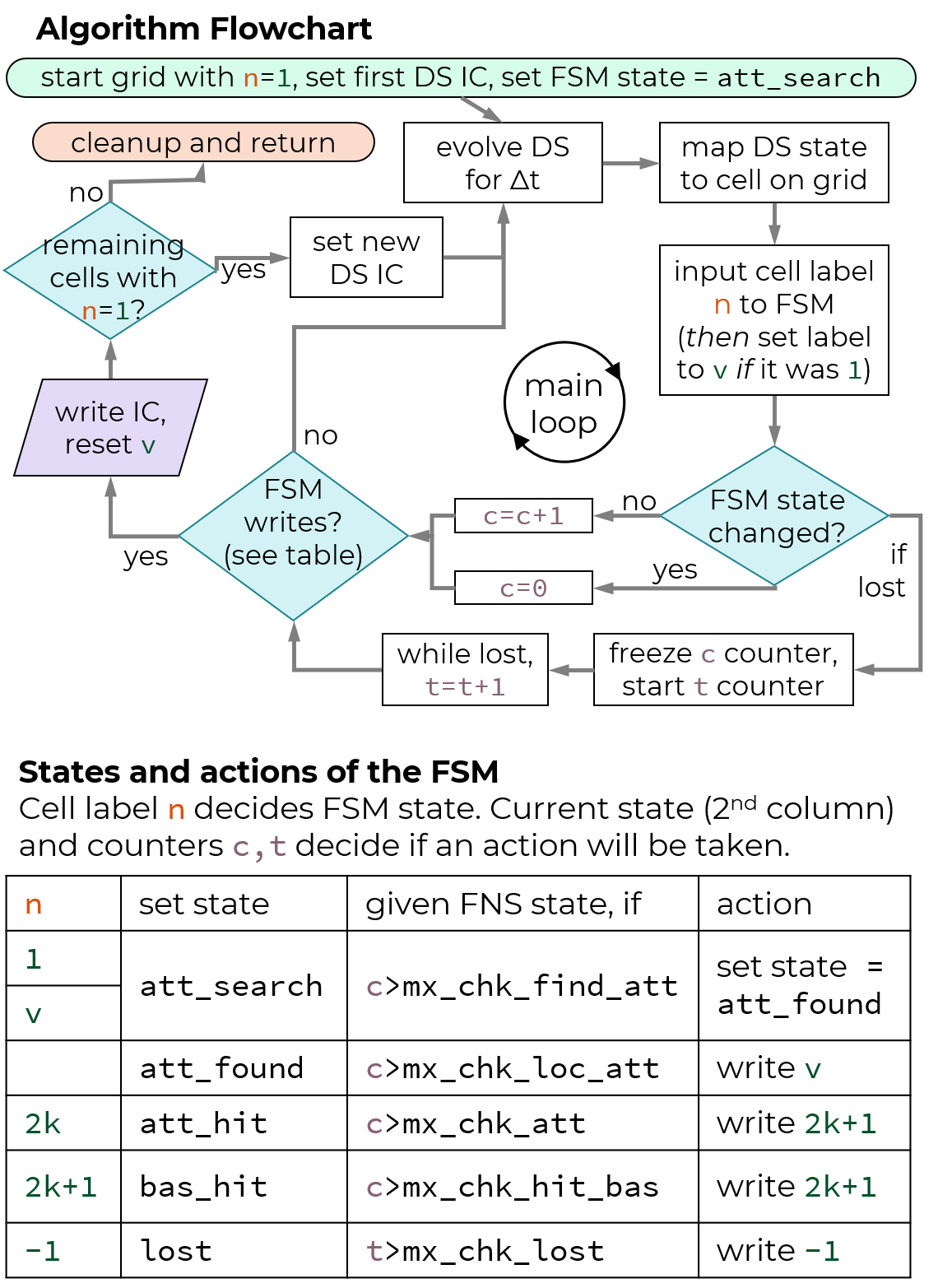}
  \caption{\label{fig:state_machine} 
  Flow diagram of our algorithm and states and actions of the finite state machine. While the FSM is on state \texttt{att\_found}, it always labels current cell as \texttt{v-1} (not shown in the flow chart). Abbreviations: DS = Dynamical System, FSM = Finite State Machine, IC = Initial Condition for DS.
}
\end{figure}

To estimate the full basins of attraction, the algorithm must identify all attractors contained in the defined grid, detect which grid cells belong to which basins, and handle the cases when a trajectory diverges or stays outside the grid. To achieve this, we propose a finite state machine (FSM) formalism built on top of the dynamical system trajectory. It coordinates the tasks of the algorithm in a systematic way and provides a flexible framework that permits new functionalities if necessary. The overall algorithm and behavior of the FSM is presented in Fig.~\ref{fig:state_machine} and in the ensuing description.

The FSM has a state and an internal counter \texttt{c}.
At each step of the main algorithm loop (Fig.~\ref{fig:state_machine}), the dynamical system is evolved for $\Delta t$ time, and its location in the state space is mapped to the enclosing grid cell. The cell label is given as the input to the FSM as shown in Fig.~\ref{fig:state_machine}.

A cell label \texttt{n} is encoded using integers. Initially every cell of the grid is labelled \texttt{1}, meaning that there is an unknown basin or attractor in this cell. The cells containing attractor points receive an even number \texttt{2k} and cells with basin points are given an odd number \texttt{2k+1} with \texttt{k>0}. Conveniently, attractors and their corresponding basins are labeled using the same \texttt{k} value, i.e., they form an even-odd pair. If the dynamical system evolution brings it outside the defined grid, \texttt{-1} is used as a cell label. Lastly, cells labelled \texttt{1} that are visited by the trajectory during the algorithm loop are labelled as \texttt{v}. We always use the next \emph{unused} odd number for \texttt{v} since it may encode the basin of attraction of a yet-to-be-identified attractor.

After receiving the cell label \texttt{n} as input, the FSM will either change its state according to the first two columns of the Table of Fig.~\ref{fig:state_machine} and set \texttt{c=0}, or stay in the same state as before and set \texttt{c=c+1}. After configuring its state and counter value, the FSM may ``write'' a value to the initial condition's cell (Fig.~\ref{fig:state_machine}), if its internal counter crosses a threshold value. After writing, the initial grid cell of the algorithm has been labelled correctly. If there are still cells labelled \texttt{1}, the process repeats with a new initial condition, otherwise the whole process terminates.

The general operation of the FSM is as follows: (1) reset counter if state/input changed, (2) increment counter while in the same state, (3) write final label to the initial grid cell once the counter is high enough (see Fig.~\ref{fig:state_machine}). This operation is independent of the actual state of the FSM. The state determines the threshold the counter must exceed for writing, and the label written in the initial cell. The default values for counter thresholds are listed in Table~\ref{tab:keywords}, while the values to be written are contained in the last column of the Table of Fig.~\ref{fig:state_machine}. The FSM has five possible states (notice that the sequence of \texttt{att\_search} $\to$ \texttt{att\_found} has been described precisely in Sect.~\ref{sec:attractor_identification}):
\begin{itemize}

\item \texttt{att\_search}: This is also the initial state of the machine and it stands for searching for an attractor.


\item \texttt{att\_found}: We have found a new attractor. This is the only state that cannot be reached via the cell label input, but rather via the state \texttt{att\_search}. In this state the FSM does not care about the input cell label. The only difference in the FSM operation is that while on state \texttt{att\_found}, the current cell is always labelled as \texttt{v-1}, which is the next unused even number, which is also the newest identified attractor. Obviously, after the end of operation of \texttt{att\_found} the numeric value for \texttt{v} is changed to \texttt{v=v+2} as we have one more new attractor in the grid.


\item \texttt{att\_hit}: The current trajectory point is in a cell containing an identified attractor point. Notice that \texttt{att\_hit} is an umbrella state: each unique attractor corresponds to a different state for the FSM. Similarly with \texttt{bas\_hit}.

\item \texttt{bas\_hit}: The input is an odd number \texttt{2k+1 < v}. Hence, the trajectory visits a cell belonging to the basin of an attractor already found. This state is not necessary for the algorithm to work but it speeds up the performance in many cases (see Sect. \ref{sec:performance}). It simply represents that if we are in the basin of attraction of a known attractor for long enough, we don't have to wait until we actually converge to the attractor to label the initial grid as belonging to the basin of said attractor.



\item \texttt{lost}: The trajectory is outside the defined grid. Here the internal counter \texttt{c} is frozen. A new counter \texttt{t} starts from \texttt{t=0} and is incremented while the FSM remains in the same state as normally. The reason for the second counter is purely for a better user experience and is not actually necessary for the algorithm to work. The second counter targets scenarios where the trajectory might slightly depart from the defined grid and return there a couple of steps later, simply because the user has not defined a large enough grid. This also means that the first counter \texttt{c} is frozen: it is not reset to 0 if the FSM returns to its prior state after after exiting the \texttt{lost} state. 
\end{itemize}

\begin{table}
\caption{Default values for the counter thresholds of each of the states of the finite state machine, see also discussion in Sect.~\ref{sec:problem_solving}.}
    \label{tab:keywords}
    \centering
    \begin{ruledtabular}
    \begin{tabular}{lr}
        Parameter & Value\\
        \hline
        \texttt{mx\_chk\_att}  & 2 \\
        \texttt{mx\_chk\_fnd\_att} & 100\\
        \texttt{mx\_chk\_loc\_att} & 100\\
        \texttt{mx\_chk\_lost} & 20 \\
        \texttt{mx\_chk\_hit\_bas}  & 10 \\
    \end{tabular}
    \end{ruledtabular}
\end{table}


The description of the algorithm above does not contain any reference to the nature of the dynamical system. The only required input is the time evolution of the state space points. As a consequence, a large variety of possible dynamical systems is admitted: discrete and continuous ones, Poincar\'e maps, and stroboscopic maps. It is also possible to track only the projected state of a dynamical system to lower-dimensional subspace of the full state space. For example, the basins of a four dimensional system can be analyzed on a projected plane of e.g., the first two variables of the system (see Sect.~\ref{sec:results} for examples). This provides a massive computational performance advantage, but is only useful in scenarios where the attractors either do not span the remaining projected dimensions, or if they do, they do not intertwine along these projected dimensions.

\subsection{Refinement of basins with known attractors}
The attractors must be contained within the limits of the defined grid when the algorithm computes their basins without prior knowledge. This is a limitation, because often one wants to focus on a region of the basins that does not contain the attractors (e.g., zooming into a part of the basins with strongly fractal boundaries as in Fig.~\ref{fig:all_basins}(c,d)). To address this, we have added a second mode of operation to the algorithm which works with user-provided already identified attractors. In this second mode, the algorithm computes the minimum distance of the current state space point versus all the attractors. When one of these distances falls bellow a given threshold $\varepsilon$, we match the initial condition with the corresponding attractor. Of course, the original algorithm can be used to first detect the attractors on a larger and coarser state space grid, which will be refined by the second mode of operation.

\subsection{Limitations and problem solving}
\label{sec:problem_solving}
Our method does not assume any approximations on the estimation of basins or attractors. In this sense it is arbitrarily precise: the more refined the grid the better the basins are estimated. Nevertheless, there are limitations and/or difficulties. The most obvious one is that localization of all attractors existing in the state space is not guaranteed for a given grid resolution. 

The total extent of the grid should be chosen large enough to actually contain the attractors fully, but also fine enough to separate attractors in state space. The step size $\Delta t$ of the integrator (in the case of continuous time systems) is also critical. It should be large enough for the trajectory to visit different cells at each step. If it is too large we may loose some performance benefits of our algorithm, but we never lose accuracy in this case. Small $\Delta t$ that make the trajectory spend several steps in the same cell in state space are a bad choice. In the code implementation we provide an automatic guess for $\Delta t$ equivalent to 10 times the average cell crossing time. 
Regarding the parameters of Table~\ref{tab:keywords}, their default values have been chosen to work well with most of the systems we tested. Obviously, increasing all of them makes the basin estimation more precise at the cost of computational performance. More specifically, these parameters should be increased in the following scenarios: \texttt{mx\_chk\_att} if attractors in the state space are very close to each other, \texttt{mx\_chk\_hit\_bas} if the basin boundaries are strongly fractal, \texttt{mx\_chk\_fnd\_att} and \texttt{mx\_chk\_loc\_att} if there are chaotic attractors.

If the algorithm does not seem to find the suspected number of attractors, or never halts because it cannot find any attractor, there are some possible actions that can help solving the problem. First increase the limits of the grid, as transients sometimes stay a long time outside the defined grid. If the dynamics is continuous, try adjusting the integrator step size and make sure the orbit visits a different cell at each step. Also the solver must be the right one for your system (e.g., stiff versus non-stiff problems).

Lastly, let us mention that finding full basins of attraction in high dimensional systems is computationally costly and strongly limited by available memory. Basin array size grows exponentially both with dimensionality and grid density, and our method needs to initialize such an array to encode the labelling described previously. Already in 10-dimensional systems with 21 grid points along each dimension this exceeds the memory available on a typical desktop computer. The best alternative we can think of is to not estimate the full basins of attraction but rather their fractions using random sampling (see discussion in end of Sect.~\ref{sec:conclusions}).

\section{Results}
\label{sec:results}
\begin{figure*}[t]
    \centering
    \includegraphics[width = \textwidth]{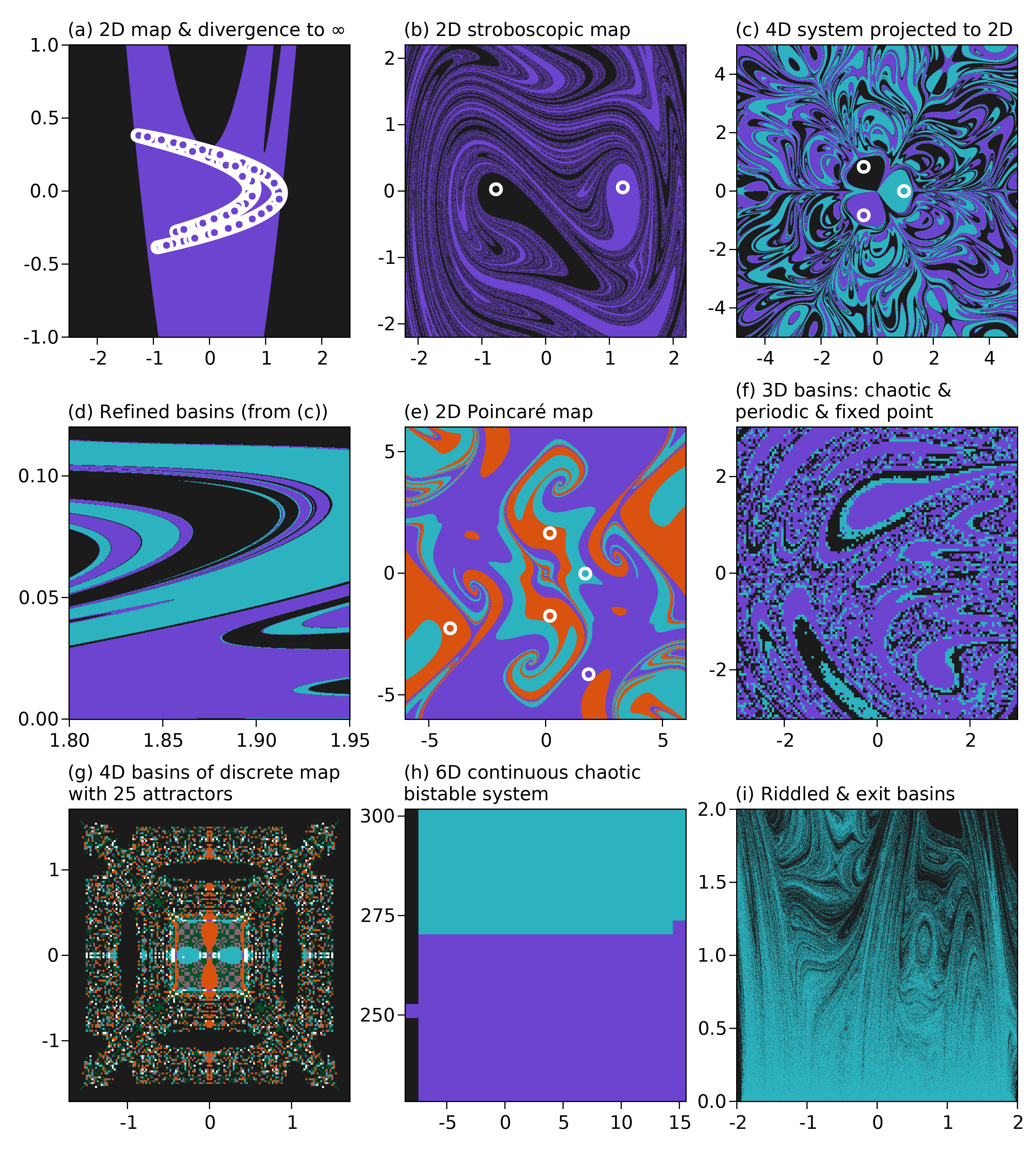}
    \caption{Basins of attraction for the scenarios discussed in Sect.~\ref{sec:results}. In the cases of (f,g,h) the basins are 3D, 4D, 6D respectively, and the plots only show a slice along two dimensions. In (a,g,h) black color corresponds to initial conditions escaping to infinity. White circles correspond to attractors for (a,b,c,e).
    In all plots the dimensions plotted are the first two of the dynamical system, except the panel (h) where it is the last two.}
    \label{fig:all_basins}
\end{figure*}

To showcase the strengths of the algorithm, we apply it to find the basins of the following scenarios:
\begin{enumerate}[label= (\alph*)]
    \item 2D discrete dynamical system with a chaotic attractor and orbits escaping to infinity (H\'enon map).
    \item 2D stroboscopic map (Duffing oscillator).
    \item 2D projection of basins of 4D continuous system with fixed points as attractors and fractal attractor basins (magnetic pendulum).
    \item Refining basins of attraction (zooming into the fractal structure) of the above.
    \item Poincar\'e map of a 3D continuous system that has interlaced attracting periodic orbits (Thomas cyclical with Poincar\'e section defined at $z=0$). On the Poincar\'e map the periodic orbits become attracting fixed points.
    \item 3D continuous dynamical system with coexistence of a chaotic attractor, an attracting periodic orbit and an attracting fixed point (Lorenz84).
    \item 4D discrete dynamical system with extreme multi-stability of $\sim 26$ coexisting attractors (nonlinearly coupled logistic maps).
    \item 6D continuous dynamical system with bistability (Lorenz96 coupled with simple energy balance model, Lorenz96EBM). One attractor is chaotic, the other periodic.
    \item 2D basins of a stroboscopic map of a forced 4D continuous system which has a basin of attraction riddled with an exit basin.
\end{enumerate}
The output is shown in Fig.~\ref{fig:all_basins}. The dynamical rule and parameters for each dynamical system is shown in Table~\ref{tab:rules}.

\begin{figure}
    \centering
    \includegraphics[width = \columnwidth]{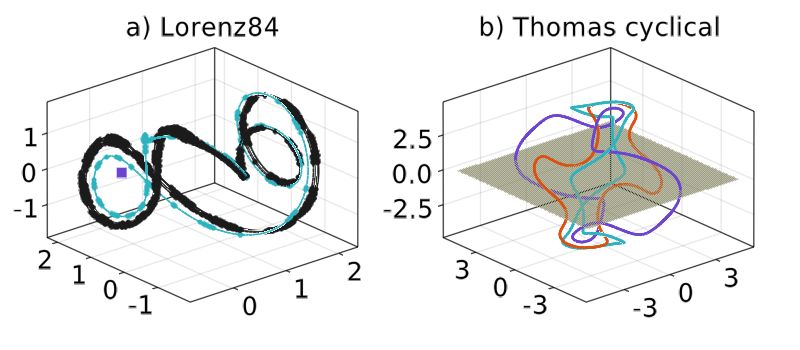}
    \caption{a) Three attractors of the Lorenz84 system (square marker for the fixed point). Circular markers are used to denote the attractor points found automatically by our algorithm, lines are used to integrate a trajectory and highlight the full attractor. b) Three periodic attractors of the Thomas cyclical (fixed point attractors also exist) and the plane used to define the Poincar\'e section.}
    \label{fig:attractors}
\end{figure}

For all cases we applied the algorithm, the expected basins have been found easily. It is especially worth it to highlight the case of Lorenz84 (Fig.~\ref{fig:all_basins}f), because two of the three attractors are extremely close in state space, see Fig.~\ref{fig:attractors}a. We used a grid of $100\times100\times100$ resolution (only a 2D slice of the full 3D basins is shown, the fraction of each basin is $\approx (0.55,0.2,0.25)$). Had we used a coarse grid resolution (less than 100 points per dimension), the two attractors would not have been separated by the algorithm. Fig.~\ref{fig:attractors}b shows the three periodic attractors of the Thomas cyclical system, and the plane used to define the Poincar\'e section. This is the plane used to produce the basins of attraction of the Poincar\'e map in Fig.~\ref{fig:all_basins}e. For the case of the 4D nonlinearly coupled logistic maps, we do not know for sure whether all attractors were found (given how many there are). There is no prior work that did a more thorough analysis on this specific system. For the 6D continuous system, the basin boundary is smooth and the two attractors are very well separated in the 6th dimension of the state space ($T$), which makes the entire process much simpler. To keep the computation time low, we used a coarse grid of $10\times10\times10\times10\times10\times101$ and only made the gridding of the last variable dense. This required about 1h30m computing time on an average machine. The basin fractions are $\approx (0.61, 0.39)$.

\section{Implementation}
\label{sec:implementation}
Our algorithm is implemented in DynamicalSystems.jl~\cite{datseris2018dynamicalsystems}. The strengths of this software, among others, are the simplicity of use and excellent numerical efficiency. Our implementation follows these principles and provides a lean interface as well as a tight computational time and memory usage. It is part of the library since v1.9. From a user perspective, using the algorithm is quite simple, and in Listing~\ref{lst:example} we present its basic application using our analysis of the Lorenz-84 model as an example.

\lstset{numbers=left,xleftmargin=2em,frame=single,framexleftmargin=1.5em}
\begin{lstlisting}[language=Python, float = *, caption = {Basic usage of our basins of attraction implementation. The listing is runnable Julia code.}, label = {lst:example}]
using DynamicalSystems, OrdinaryDiffEq
# Create instance of `DynamicalSystem`:
function lorenz84(u, p, t)
    F, G, a, b = p
    x, y, z = u
    dx = -y^2 -z^2 -a*x + a*F
    dy = x*y - y - b*x*z + G
    dz = b*x*y + x*z - z
    return SVector(dx, dy, dz)
end
F, G, a, b = 6.886, 1.347, 0.255, 4.0
p  = [F, G, a, b]
u0 = rand(3) # initial condition doesn't matter
lo = ContinuousDynamicalSystem(lorenz84, u0, p)
# Calculate basins of attraction
xg = range(-1, 3;   length = 100)
yg = range(-2, 3;   length = 100)
zg = range(-2, 2.5; length = 100)
grid = (xg, yg, zg)
diffeq = (alg = Vern9(), reltol = 1e-9, abstol = 1e-9)
basins, attractors = basins_of_attraction(grid, lo; diffeq)
# Further use output for e.g., Lyapunov exponents or basin fractions:
fracs = basin_fractions(basins)
for (key, att) in attractors
    u0 = att[1] # First found point of attractor
    ls = lyapunovspectrum(lo, 10000; u0)
    println("Attractor $(key) has spectrum: $(ls) and fraction: $(fracs[key])")
end

\end{lstlisting}

The user first needs to define a \verb#DynamicalSystem# instance, done in lines 3-14 of the listing. Then, with the appropriate grid of initial conditions, the function \verb#basins_of_attraction# is called as listed in lines 16-21. The first output of the function is an array \verb#basins# with size identical to the grid. Its elements are the IDs of the attractor labelling each initial condition. The second output \verb#attractors# is a dictionary, mapping attractor IDs to the automatically estimated attractor points in the state space. These points have the dimensionality of the state space which could be higher than that of the grid. The function \verb#basins_of_attraction# allows for several keywords including those of Table~\ref{tab:keywords}.

\subsection{Integration with DynamicalSystems.jl and the Julia ecosystem}
Implementing our algorithm in DynamicalSystems.jl instead of an isolated piece of software comes with big advantages, the first being the simplicity and high-levelness of Listing~\ref{lst:example}. More importantly though, our implementation is able to communicate and be used with the rest of the library, and in fact the whole Julia ecosystem, directly. For example, in lines 24-28 of the Listing we \emph{reuse} the existing defined structures \verb#lo# and \verb#attractors# to calculate the Lyapunov exponents of each attractor. The output \verb#basins# can be further used with functions of the library such as \verb#basin_fractions#, \verb#tipping_probabilities# or \verb#basin_entropy#. These measures are useful in the analysis of dynamical systems in terms of basin stability~\cite{menck2013basin}, tipping probabilities~\cite{kaszas2019tipping} or basin entropy\cite{daza2016basin}. Lastly, DynamicalSystems.jl integrates with the Julia library DifferentialEquations.jl~\cite{DifferentialEquations.jl}. Users can pick any of the hundreds of ODE solvers from this library and adjust on the fly any accuracy-related option by providing the extra keyword \verb#diffeq#. In our work we used Verner's ``Most Efficient'' 9/8 Runge-Kutta solver with strict error tolerances by providing the keyword \verb#diffeq# as shown in line 20 of the Listing.

\subsection{Performance considerations}
\label{sec:performance}
Julia, its suite of differential equations solvers, and the optimizations of DynamicalSystems.jl, provide excellent numeric performance that our implementation takes advantage of. For example, estimating the 3D basins of attraction of the Lorenz-84 system for a $80\times 80\times 80$ grid resolution (512000 initial conditions)  requires 3 minutes on a medium performance computer with CPU AMD Ryzen 5 3600 6-Core (only one core is used as our method is not parallelizable).

\begin{figure}
    \centering
    \includegraphics[width = \columnwidth]{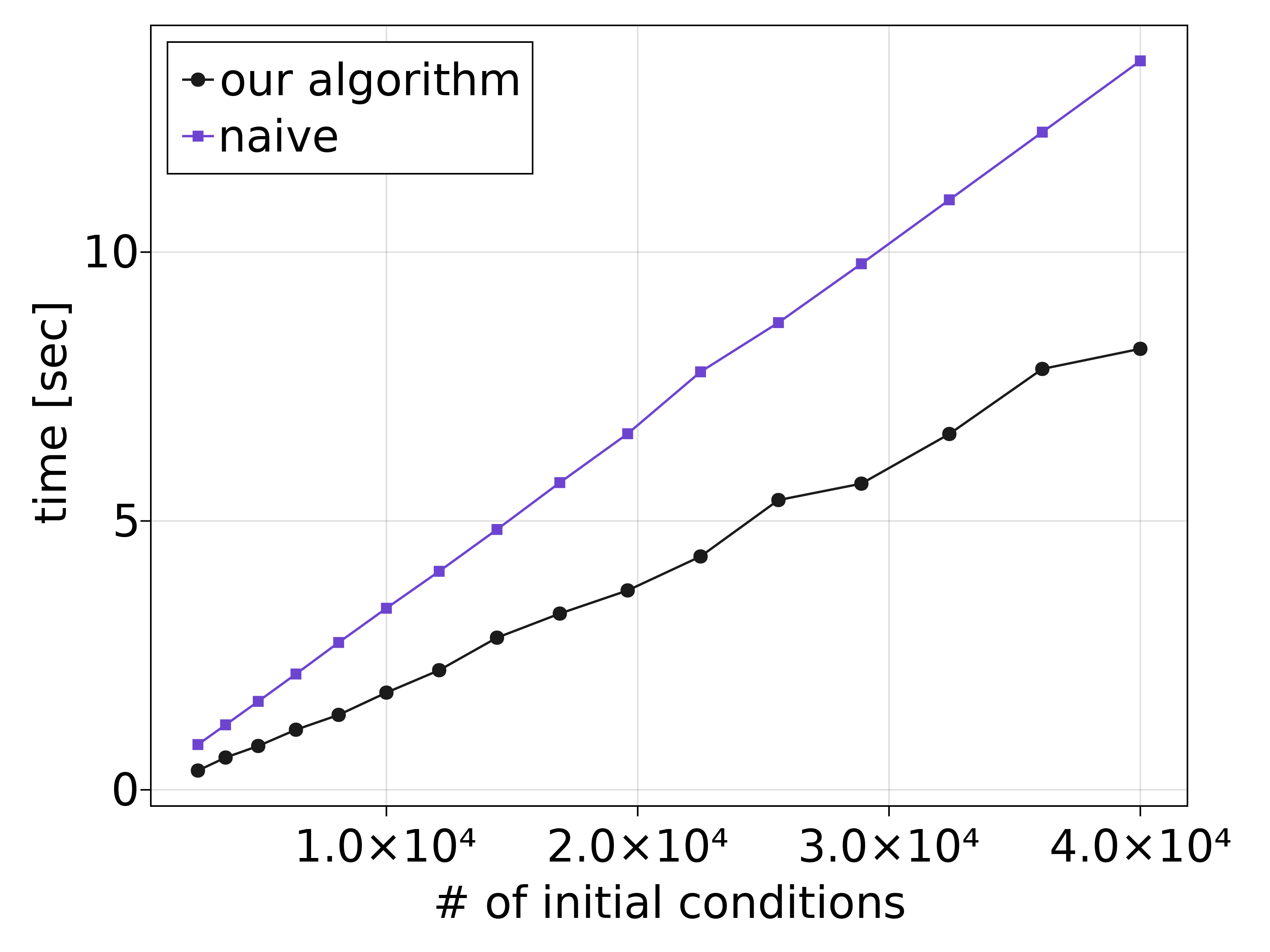}
    \caption{Benchmark comparison of creating Fig.~\ref{fig:all_basins}(c) using our algorithm or the naive approach. The timings of the latter do not include any consideration of the time needed to identify and catalogue the attractors while this is included in our algorithm.}
    \label{fig:benchmark}
\end{figure}

To obtain a language-agnostic performance estimate of our algorithm, we will compare the computation of Fig.~\ref{fig:all_basins}(c) using our algorithm against the naive approach where each initial condition is integrated until convergence to a fixed point and later mapped to one of the known three attracting fixed points. The case of Fig.~\ref{fig:all_basins}(c) is, by choice, the most unfair case we could have chosen for such a comparison: (i) the attractors are fixed points, the easiest (and perhaps only) kind of attractors the naive approach can find, (ii) the basins of attraction are strongly fractal, which reduces some of our algorithm's performance benefits. Nevertheless, as we can see in Fig.~\ref{fig:benchmark}, our method outperforms the naive approach even when excluding any time necessary to actually find the attractors (which could well be the hardest step depending on the occasion).

One of the reasons for this improvement is the use of the information already stored in the grid. The algorithm checks if the trajectory visits cells labelled as basins. If it is the case for \texttt{mx\_chk\_hit\_bas} times in a row for the same basin, the initial condition belongs to this basin. As the grid is filled, this event becomes more and more frequent and shortens the time of identification.

\section{Conclusions}
\label{sec:conclusions}

The automatic estimation of attractors and their basins of attraction is not an easy task for nonlinear dynamical systems. In this work we presented an algorithm that does better than previous solutions. It is based on a definition of an appropriate finite state machine on the state space, whose desired operation is guaranteed by the Poincar\'e recurrence theorem. The algorithm is straightforward to use, computationally performant and is implemented in the general purpose library DynamicalSystems.jl. In Sect.~\ref{sec:results} we applied our algorithm to a large array of different scenarios and demonstrated its success with all of them. We cannot underestimate the importance of numerical methods in the field of nonlinear dynamics. Our paper aims at completing the basic toolbox of researchers with a ready-to-use and versatile tool for estimating attractors and their basins of attraction.

In the near future we will enrich this functionality with a recent approach for the estimation of the basin fractions without computing the full basins of attraction from Stender \& Hoffmann~\cite{stender2021bstab}, called bSTAB. This method transforms a trajectory into a vector of features, for example the mean and the variance of the timeseries, for its later classification in the feature space. It is an interesting technique that does not require a huge in-memory matrix initialization, but it requires the user to have a basic idea of the attractors already, as well as which features can be used to distinguish between them. We plan to implement this method in DynamicalSystems.jl soon, leveraging our existing algorithm to estimate the basins of attraction.

\section*{Acknowledgments}
A.W. acknowledges the support  from the Spanish State Research Agency (AEI) and the European Regional Development Fund (ERDF, EU) under project PID2019-105554GB-I00.
\section*{Code Availability}
Besides the open source implementation, the entire code base required to produce the figures for this work is also available as open source code online\cite{Codebase}.
\section*{Data Availability}
This work does not have associated data besides the output produced from the provided code base.

\renewcommand{\arraystretch}{2.0}
\begin{table*}[t]
  \caption{Dynamical rules and parameters for systems used. For the magnetic pendulum the magnet locations $\mathbf{x}_i$ are equispaced on the unit circle. For the coupled logistic maps $u_i'$ denotes the next state and $i$ runs from 1 to $D$ (the state space dimensionality).}
  \label{tab:rules}
  \centering
  \begin{tabular*}{\textwidth}{l @{\extracolsep{\fill}}  l @{\extracolsep{\fill}}  l}
  \toprule
  System & Dynamical rule & Parameters \\ \hline
    H\'enon map\cite{Henon1976} & $ x_{n+1} = 1 - ax^2_n+y_n,\quad y_{n+1} = bx_n$ & $a=1.4, b=0.3$ \\ \hline
    Duffing oscillator\cite{Kanamaru2008} & $\ddot{x} + d \dot{x} + \beta x + x^3 = f \cos(\omega t)$ & $\omega = 1.0, f = 0.2, d = 0.15, \beta = -1.0$ \\ \hline
    Magnetic pendulum &  $\begin{aligned}
\ddot{\mathbf{x}} &= -\omega ^2\mathbf{x} - \alpha \dot{\mathbf{x}} - \sum_{i=1}^N \frac{ \mathbf{x} - \mathbf{x}_i}{D_i^3},\quad \mathbf{x} = (x,y) \\
D_i &= \sqrt{(x-x_i)^2  + (y-y_i)^2 + d^2}
\end{aligned}$ & $\alpha = 0.2, \omega = 1, d = 0.3, N=3$ \\ \hline
Thomas cyclical\cite{THOMAS1999} & $\dot{x} = \sin(y) - bx, \quad \dot{y} = \sin(z) - by,\quad \dot{z} = \sin(x) - bz$ & $b = 0.1665$ \\ \hline
Lorenz84\cite{Freire2008} & $
\dot x = - y^2 - z^2 - ax + aF, \, \dot y = xy - y - bxz + G, \, \dot z = bxy + xz - z$ & $F = 6.886, G = 1.347, a = 0.255, b = 4.0$\\ \hline
Coupled logistic maps\cite{Bezruchko2003} & $u_i' = \lambda - u_i^2 + k \sum_{j\ne i}^D (u_j^2 - u_i^2)$ & $D = 4, \lambda = 1.2, k = 0.08$ \\ \hline
Lorenz96EBM\cite{Gelbrecht2021} &  $\begin{aligned}
\dot{x}_i &= (x_{i+1} - x_{i-2})x_{i-1} - x_i + F\left(1+\beta\frac{T-\bar{T}}{\Delta_T}\right) \\
\dot{T} &= S\left(1 - a_0 + 0.5a_1\tanh(T-\bar{T}) \right) - \sigma T^4 - \alpha \left(\frac{\mathcal{E}(\mathbf{x})}{0.6F^{1.33} - 1} \right) \\
\mathcal{E}(\mathbf{x}) & = \frac{1}{2N}\sum_{i=1}^N x_i^2
\end{aligned}$ & $\begin{aligned}
N &= 5, F = 8, S = 8, a_0 = 0.5\\
a_1 &= 0.4, \bar{T} = 270, \Delta_T = 60, \alpha = 2,\\
\beta &= 1, \sigma = 1/180
\end{aligned}$ \\ \hline
Riddled system\cite{Ott1994Riddled} &  $\begin{aligned}
\dot{x} &= v_x, \quad \dot{y} = v_z \\
\dot{v}_x &= -\gamma v_x - ( -4x(1-x^2) +y^2) + f_0 \sin(\omega t)x_0 \\
\dot{v}_y &= -\gamma v_y - (2y(x+\bar{x})) + f_0 \sin(\omega t)y_0
\end{aligned}$ & $\begin{aligned}
\gamma &= 0.05, \bar{x} = 1.9, f_0 = 2.3\\
\omega &= 3.5, x_0 = 1, y_0 = 0
\end{aligned}$ \\
\toprule
  \end{tabular*}
\end{table*}

\nocite{*}
\bibliography{REFERENCES}

\end{document}